\documentclass{amsart}

\usepackage{amssymb}
\usepackage{mathtools}
\usepackage{mathrsfs}

\usepackage{enumitem}
\usepackage{diagbox}

\usepackage[colorlinks=true,
            allcolors=black]{hyperref}

\newtheorem{theorem}{Theorem}[section]
\newtheorem{lemma}[theorem]{Lemma}
\newtheorem{proposition}[theorem]{Proposition}
\newtheorem{corollary}[theorem]{Corollary}

\theoremstyle{definition}

\newtheorem{example}[theorem]{Example}

\numberwithin{equation}{section}

\newcommand{\norm}[1]{|#1| \mkern 0.5mu}
\newcommand{\ind}{\mathrel{\perp\!\!\!\perp}}

\title[Maximum distance covariance]{The maximum distance covariance of bounded random vectors}
\author{John \c{C}amk{\i}ran}
\address{University of Toronto, Toronto, Ontario, M5S 3E4, Canada}
\email{john.camkiran@utoronto.ca}

\begin{document}

\begin{abstract}
A sharp upper bound is obtained for the distance covariance of bounded random vectors in terms of the diameters of their supports alone. The bound is attained for all pairs of support diameter constraints and hence provides the maximum distance covariance subject to constraints of this kind. A normalized measure of dependence is derived by dividing distance covariance by the sharp bound when it is positive. Several properties of the resulting quotient are proved, including domination by distance correlation, invariance under similarities, and subadditivity under sums of independent pairs. The quotient is shown to be nonincreasing under the addition of independent noise, a property not shared by distance correlation. The usual sample distance covariance inherits the sharp bound from the population quantity and the associated normalized statistic is strongly consistent.
\end{abstract}

\maketitle

\section{Introduction}
\label{sec:introduction}

Dependence measures are an important object of study in mathematical statistics \cite{Renyi1959, Mori2019}. The classical measure of dependence between two random variables $X$ and $Y$, the product-moment covariance $\operatorname{cov}(X,Y)$, has the limitation that $\operatorname{cov}(X,Y) = 0$ does not characterize independence. To give a standard counterexample, if $X$ is symmetric about zero with $\mathbb{E}\lvert X\rvert^3 < \infty$ and $Y = X^2$ is nondegenerate, then $\mathbb{E}[X] = \mathbb{E}[XY] = 0$ and so $\operatorname{cov}(X,Y) = 0$ despite $Y$ being a function of $X$.

This problem was addressed by Sz\'{e}kely, Rizzo, and Bakirov with the introduction of \emph{distance covariance} $\operatorname{dCov}(X,Y)$, which is zero if and only if $X$ and $Y$ are independent \cite{Szekely2007}. Distance covariance is an energy statistic \cite{Szekely2013b, Rizzo2016} originally defined for random vectors $X$ in $\mathbb{R}^N$ and $Y$ in $\mathbb{R}^M$ with finite first moments as the nonnegative square root of
\[
    \operatorname{dCov}^2(X, Y) = \|f_{X, Y}(s, t) - f_X(s) f_Y(t)\|^2_w,
\]
where $f_X$ and $f_Y$ are the marginal characteristic functions, $f_{X, Y}$ is the joint characteristic function, and $\|\cdot\|_w$ is the $L^2$-norm induced by the weight function proposed in \cite{Szekely2007} and characterized in \cite{Szekely2012}. 

Distance covariance, which coincides with Brownian covariance when the second moments of $X$ and $Y$ are also finite \cite{Szekely2009}, satisfies the following Cauchy--Schwarz-type inequality:
\begin{equation}
    \operatorname{dCov}(X,Y) \leq \sqrt{\operatorname{dCov}(X,X)\operatorname{dCov}(Y,Y)}.
    \label{ineq:cauchy-schwarz}
\end{equation}
By analogy with Pearson's product-moment correlation, Sz\'{e}kely, Rizzo, and Bakirov defined a normalization of distance covariance termed \emph{distance correlation} given by
\begin{equation}
    \operatorname{dCor}(X,Y) = \frac{\operatorname{dCov}(X,Y)}{\sqrt{\operatorname{dCov}(X,X)\operatorname{dCov}(Y,Y)}}
    \label{eq:distance-correlation}
\end{equation}
whenever $\operatorname{dCov}(X,X) \operatorname{dCov}(Y,Y) > 0$ and $\operatorname{dCor}(X,Y) = 0$ otherwise \cite{Szekely2007}. It may be noted that both of the foregoing relations involve the marginal distributions of $X$ and $Y$.

It is natural to ask what can be said about $\operatorname{dCov}(X,Y)$ when only rough support information is available, such as the diameters of the supports of $X$ and $Y$. Let $\mathcal{A},\mathcal{B}\geq 0$ be two real numbers satisfying $\operatorname{diam}(\operatorname{supp} X)\leq \mathcal{A}$ and $\operatorname{diam}(\operatorname{supp} Y)\leq \mathcal{B}$. Then the distance of each random vector to an independent copy of itself is uniformly bounded from above by $\mathcal{A}$ and $\mathcal{B}$ respectively. It follows that the double-centered versions of these distances, which enter the energy-statistical representation of $\operatorname{dCov}^2(X,Y)$ \cite{Szekely2013b}, are bounded in absolute value by $2\mathcal{A}$ and $2\mathcal{B}$ respectively. One thus obtains the crude diameter-only estimate $\operatorname{dCov}(X,Y)\leq 2\sqrt{\mathcal{A}\mathcal{B}}$. The second inequality in Proposition~2.3 of Lyons\footnote{Lyons writes $\operatorname{dcov}(X,Y)$ for the quantity here denoted by $\operatorname{dCov}^2(X,Y)$.} \cite{Lyons2013} improves this estimate to $\operatorname{dCov}(X,Y)\leq \sqrt{\mathcal{A}\mathcal{B}}$.

In this paper we establish a sharp upper bound on distance covariance that, much like the bounds in the preceding paragraph, depends only on the diameters of the supports of the underlying random vectors. Unlike those estimates, ours is attained by some pair $(X,Y)$ of random vectors for every pair $(\mathcal{A}, \mathcal{B})$ of support-diameter bounds and thus gives the maximum distance covariance under the corresponding constraints. The following theorem is proved in Section~\ref{sec:proof-of-main-theorem}.

\begin{theorem}
\label{thm:main}
For $\mathcal{A}, \mathcal{B} \geq 0$ let $\mathscr{C}(\mathcal{A},\mathcal{B})$ be the class of random vector pairs $(X, Y)$ such that
\[
    \operatorname{diam}(\operatorname{supp} X) \leq \mathcal{A}, \qquad
    \operatorname{diam}(\operatorname{supp} Y) \leq \mathcal{B}.
\]
Then
\[
    \max_{(X,Y) \in \mathscr{C}(\mathcal{A},\mathcal{B})} \operatorname{dCov}(X, Y) = \frac{1}{2} \sqrt{\mathcal{A} \mathcal{B}}.
\]
\end{theorem}

In Section~\ref{sec:normalized-dependence-measure} we use our bound to define a normalized measure of dependence that combines the scale invariance of distance correlation with the subadditivity of distance covariance under sums of independent pairs. The latter property has the consequence that the normalized measure is nonincreasing under the addition of independent noise, which turns out not to be the case for distance correlation. In Section~\ref{sec:sample-bound-normalization} we show that the sharp bound and the related normalization carry over to the usual V-statistic-based sample distance covariance \cite{Szekely2023} and briefly consider the behavior of the unbiased U-statistic for the squared distance covariance.

Before proceeding with the paper we introduce some notation and conventions. The symbols $X$ and $Y$ are used to denote bounded random vectors taking values in Euclidean spaces of possibly different dimensions. The pairs $(X,Y)$, $(X',Y')$, and $(X'',Y'')$ are mutually independent and identically distributed (i.i.d.). Finally we write $\norm{\cdot}$ for the Euclidean norm, $\operatorname{supp}(\cdot)$ for the support of (the law of) a random vector, and $\operatorname{diam}(\cdot)$ for the diameter of a bounded set in Euclidean space.

\section{Proof of Theorem~\ref{thm:main}}
\label{sec:proof-of-main-theorem}

From Remark~3 of \cite{Szekely2007} we have the identity
\begin{multline}
\label{eq:dcov-identity}
    \operatorname{dCov}^2(X, Y)
    = \mathbb{E}[\norm{X-X'}\norm{Y-Y'}]
    +\mathbb{E}[\norm{X-X'}] \mathbb{E}[\norm{Y-Y'}] \\
    - 2\mathbb{E}[\norm{X-X'}\norm{Y-Y''}].
\end{multline}
Our first lemma uses \eqref{eq:dcov-identity} to establish the diameter-only upper bound.

\begin{lemma}
\label{lem:sharp-bound}
Let $\mathcal{A},\mathcal{B}\geq 0$. If $\operatorname{diam}(\operatorname{supp} X) \leq \mathcal{A}$ and $\operatorname{diam}(\operatorname{supp} Y) \leq \mathcal{B}$ then
\[
    \operatorname{dCov}(X, Y) \leq \frac{1}{2} \sqrt{\mathcal{A}\mathcal{B}}.
\]
\end{lemma}

\begin{proof}
Set
\[
    D_1 := \norm{X-X'}, \qquad D_2 := \norm{X-X''}.
\]
Since $X'$ and $X''$ are independent copies of $X$, we have
\[
    \mathbb{E}[D_1D_2] = \mathbb{E}[\mathbb{E}[D_1D_2\mid X]] = \mathbb{E}[\mathbb{E}[D_1\mid X]\mathbb{E}[D_2\mid X]] = \mathbb{E}[\mathbb{E}[D_1\mid X]^2].
\]
Subtracting $\mathbb{E}[D_1]^2 = \mathbb{E}[\mathbb{E}[D_1\mid X]]^2$ from the left- and right-hand sides yields
\[
    \mathbb{E}[D_1D_2] - \mathbb{E}[D_1]^2
    =
    \mathbb{E}[\mathbb{E}[D_1\mid X]^2]
    -\mathbb{E}[\mathbb{E}[D_1\mid X]]^2
    =
    \operatorname{var}(\mathbb{E}[D_1\mid X])
    \geq 0.
\]
Hence $\mathbb{E}[D_1D_2] \geq \mathbb{E}[D_1]^2$. Applying \eqref{eq:dcov-identity} with $Y=X$ we get
\begin{equation}
\label{ineq:dcov-variance-bound}
    \operatorname{dCov}^2(X,X) = \mathbb{E}[D_1^2] + \mathbb{E}[D_1]^2 - 2\mathbb{E}[D_1D_2] \leq \mathbb{E}[D_1^2] - \mathbb{E}[D_1]^2 = \operatorname{var}(D_1).
\end{equation}
It follows from $\operatorname{diam}(\operatorname{supp} X) \leq \mathcal{A}$ that $0 \leq D_1=\norm{X-X'} \leq \mathcal{A}$ a.s. By Popoviciu's inequality \cite{Popoviciu1935}
\begin{equation}
\label{ineq:var-bound}
    \operatorname{var}(D_1) \leq \frac{\mathcal{A}^2}{4}.
\end{equation}
Combining \eqref{ineq:dcov-variance-bound} and \eqref{ineq:var-bound} and arguing similarly with $Y$ in place of $X$ gives
\begin{equation}
\label{ineq:self-bounds}
    \operatorname{dCov}(X,X) \leq \frac{\mathcal{A}}{2},
    \qquad
    \operatorname{dCov}(Y,Y) \leq \frac{\mathcal{B}}{2}.
\end{equation}
Using \eqref{ineq:cauchy-schwarz} and \eqref{ineq:self-bounds} we obtain
\[
    \operatorname{dCov}(X, Y) \leq \sqrt{\operatorname{dCov}(X,X) \operatorname{dCov}(Y,Y)} \leq \frac{1}{2} \sqrt{\mathcal{A} \mathcal{B}}.
\]
\end{proof}

Inspecting the proof of Lemma~\ref{lem:sharp-bound}, one observes that for $\mathcal{A}, \mathcal{B}>0$, the stated bound is attained if and only if the following equalities hold:
\[
    \operatorname{dCov}(X,X)=\frac{\mathcal{A}}{2},
    \qquad
    \operatorname{dCov}(Y,Y)=\frac{\mathcal{B}}{2},
    \qquad
    \operatorname{dCor}(X,Y)=1.
\]
To prove the next lemma, we construct a simple two-point distribution that meets these conditions for any compact sets $\mathcal{X}, \mathcal{Y}$ with positive diameters.

\begin{lemma}
\label{lem:attainment}
For any two nonempty compact sets $\mathcal{X} \subseteq \mathbb{R}^N$ and $\mathcal{Y} \subseteq \mathbb{R}^M$ such that $\operatorname{diam}(\mathcal{X}) = \mathcal{A}$ and $\operatorname{diam}(\mathcal{Y}) = \mathcal{B}$, there exist random vectors $X$ and $Y$ with $\operatorname{supp} X \subseteq \mathcal{X}$ and $\operatorname{supp} Y \subseteq \mathcal{Y}$ satisfying
\[
    \operatorname{dCov}(X, Y) = \frac{1}{2}\sqrt{\mathcal{A}\mathcal{B}}.
\]
\end{lemma}

\begin{proof}
Since $\mathcal{X}$ and $\mathcal{Y}$ are compact the continuous functions
\[
(x,x') \mapsto \norm{x-x'}, \qquad (y,y') \mapsto \norm{y-y'}
\]
attain their maxima on $\mathcal{X} \times \mathcal{X}$ and $\mathcal{Y} \times \mathcal{Y}$ respectively. In other words there exist $x_1,x_2 \in \mathcal{X}$ and $y_1,y_2 \in \mathcal{Y}$ such that
\[
    \norm{x_1 - x_2} = \mathcal{A}, \qquad \norm{y_1 - y_2} = \mathcal{B}.
\]

If $\mathcal{A}\mathcal{B} = 0$ the claim is trivial. Assume $\mathcal{A}\mathcal{B} > 0$. We show that the bound is attained for random vectors $X$ and $Y$ of the forms
\begin{equation}
\label{eq:two-point-construction}
    X = 
    \begin{cases}
        x_1 & \text{w.p. } 1/2 \\
        x_2 & \text{w.p. } 1/2,
    \end{cases}
    \qquad
    Y = 
    \begin{cases}
        y_1 & \text{if } X = x_1 \\
        y_2 & \text{if } X = x_2,
    \end{cases}
\end{equation}
in which case $(X,Y)$ takes each of the values $(x_1,y_1)$ and $(x_2,y_2)$ with probability $1/2$. 

By \eqref{eq:two-point-construction} the products $\norm{X-X'}\norm{Y-Y'}$ and $\norm{X-X'}\norm{Y-Y''}$ equal $\mathcal{A}\mathcal{B}$ on the events $\{X\neq X'\}$ and $\{X \neq X', X \neq X''\}$ respectively and vanish otherwise. Since $X,X',X''$ are independent and each takes its two values with probability $1/2$, we have
\[
    \Pr(X\neq X') = \frac{1}{2},
    \qquad
    \Pr(X\neq X', X\neq X'') = \frac{1}{4}.
\]
Consequently
\begin{align*}
    \mathbb{E}[\norm{X-X'}\norm{Y-Y'}]
    &= \frac{1}{2}\mathcal{A}\mathcal{B},\\
    \mathbb{E}[\norm{X-X'}]\mathbb{E}[\norm{Y-Y'}]
    &= \frac{1}{4}\mathcal{A}\mathcal{B},\\
    \mathbb{E}[\norm{X-X'}\norm{Y-Y''}]
    &= \frac{1}{4}\mathcal{A}\mathcal{B}.
\end{align*}

Substituting these values into \eqref{eq:dcov-identity} we obtain
\[
    \operatorname{dCov}^2(X,Y) = \frac{1}{2} \mathcal{A}\mathcal{B} + \frac{1}{4} \mathcal{A}\mathcal{B} - \frac{2}{4} \mathcal{A}\mathcal{B} = \frac{1}{4} \mathcal{A}\mathcal{B}.
\]
Therefore
\[
    \operatorname{dCov}(X,Y) = \frac{1}{2} \sqrt{\mathcal{A}\mathcal{B}}
\]
as claimed.
\end{proof}

\begin{proof}[Proof of Theorem~\ref{thm:main}]
Lemma~\ref{lem:sharp-bound} asserts that
\[
    \operatorname{dCov}(X,Y)\leq \frac{1}{2}\sqrt{\mathcal{A}\mathcal{B}}
\]
for every $(X,Y)\in\mathscr{C}(\mathcal{A},\mathcal{B})$, and hence $\frac{1}{2}\sqrt{\mathcal{A}\mathcal{B}}$ is an upper bound for $\operatorname{dCov}$ on $\mathscr{C}(\mathcal{A},\mathcal{B})$. Choose compact sets $\widetilde{\mathcal{X}}$ and $\widetilde{\mathcal{Y}}$ with $\operatorname{diam}(\widetilde{\mathcal{X}})=\mathcal{A}$ and $\operatorname{diam}(\widetilde{\mathcal{Y}})=\mathcal{B}$ (one can always take $\widetilde{\mathcal{X}} = [0,\mathcal{A}]$ and $\widetilde{\mathcal{Y}} = [0,\mathcal{B}]$). Lemma~\ref{lem:attainment} provides random vectors $\widetilde{X}$ supported in $\widetilde{\mathcal{X}}$ and $\widetilde{Y}$ supported in $\widetilde{\mathcal{Y}}$ satisfying $\operatorname{dCov}(\widetilde{X},\widetilde{Y})=\frac{1}{2}\sqrt{\mathcal{A}\mathcal{B}}$. Since $(\widetilde{X},\widetilde{Y})\in\mathscr{C}(\mathcal{A},\mathcal{B})$, the upper bound is attained and $\frac{1}{2}\sqrt{\mathcal{A}\mathcal{B}}$ is the maximum of $\operatorname{dCov}$ on $\mathscr{C}(\mathcal{A},\mathcal{B})$.
\end{proof}

\section{Normalized measure of dependence}
\label{sec:normalized-dependence-measure}

A normalized measure of dependence is obtained by dividing distance covariance by the sharp bound of Lemma~\ref{lem:sharp-bound} when it is positive. Set
\[
    \mathcal{A}:=\operatorname{diam}(\operatorname{supp} X),
    \qquad
    \mathcal{B}:=\operatorname{diam}(\operatorname{supp} Y).
\]
Define
\begin{equation}
\label{eq:J-definition}
    J(X, Y) = 
    \begin{cases}
        \dfrac{2\operatorname{dCov}(X, Y)}{\sqrt{\mathcal{A}\mathcal{B}}} & \mathcal{A}\mathcal{B}>0 \\
        0 & \text{otherwise}.
    \end{cases}
\end{equation}

We establish several properties of $J$. In the following proposition, $Z$ and $W$ denote bounded random vectors of the same dimensions as $X$ and $Y$ respectively.
\begin{proposition}\leavevmode
\label{prop:J-properties}
\begin{enumerate}[label=(\roman*)]
    \item $0 \leq J(X,Y) \leq 1$.
    \item $J(X,Y) = J(Y,X)$.
    \item $J(X,Y) = 0 \Longleftrightarrow X \ind Y$.
    \item $J(X,Y) \leq \operatorname{dCor}(X,Y)$.
    \item $J(X,Y) = J(T(X), S(Y))$ for all similarities $T$ and $S$.
    \item $J(X+Z,Y+W) \leq J(X,Y)+J(Z,W)$ if $(X,Y) \ind (Z,W)$.
\end{enumerate}
\end{proposition}

\begin{proof}
(i) is immediate from Lemma~\ref{lem:sharp-bound}; (ii) is a consequence of the symmetry $\operatorname{dCov}(X,Y)=\operatorname{dCov}(Y,X)$; (iii) follows from $J(X,Y) = 0 \Longleftrightarrow \operatorname{dCov}(X,Y)=0$ and $\operatorname{dCov}(X,Y)=0 \Longleftrightarrow X\ind Y$ \cite{Szekely2007}. It remains to show (iv)--(vi). 

\bigskip
\noindent \textbf{(iv).} If $\operatorname{dCov}(X,Y) = 0$ then $J(X,Y) = 0 \leq \operatorname{dCor}(X,Y)$. Assume $\operatorname{dCov}(X,Y) \mkern-1.5mu > \mkern-1.5mu 0$, which implies that $\operatorname{dCov}(X,X) \operatorname{dCov}(Y,Y) > 0$ and $\mathcal{A}\mathcal{B}>0$. Lemma~\ref{lem:sharp-bound} applied to $(X,X)$ and $(Y,Y)$ yields
\[
    \operatorname{dCov}(X,X)\leq \frac{\mathcal{A}}{2},
    \qquad
    \operatorname{dCov}(Y,Y)\leq \frac{\mathcal{B}}{2},
\]
and hence
\[
    \frac{2\sqrt{\operatorname{dCov}(X,X)\operatorname{dCov}(Y,Y)}}{\sqrt{\mathcal{A}\mathcal{B}}}\leq 1.
\]
Using \eqref{eq:distance-correlation} we compute
\begin{align*}
    \frac{J(X,Y)}{\operatorname{dCor}(X,Y)}
    &= \frac{2\operatorname{dCov}(X,Y)}{\sqrt{\mathcal{A}\mathcal{B}}} \frac{\sqrt{\operatorname{dCov}(X,X)\operatorname{dCov}(Y,Y)}}{\operatorname{dCov}(X,Y)} \\
    &= \frac{2\sqrt{\operatorname{dCov}(X,X)\operatorname{dCov}(Y,Y)}}{\sqrt{\mathcal{A}\mathcal{B}}} \leq 1.
\end{align*}
Thus $J(X,Y)\leq \operatorname{dCor}(X,Y)$.

\bigskip
\noindent \textbf{(v).} Write $T(x) = \alpha Qx + a$ and $S(y) = \beta Ry + b$ for positive scalars $\alpha, \beta > 0$, vectors $a \in \mathbb{R}^N$ and $b \in \mathbb{R}^M$, and orthogonal matrices $Q \in \mathbb{R}^{N \times N}$ and $R \in \mathbb{R}^{M \times M}$. If $\mathcal{A}\mathcal{B}=0$ then $J(X,Y) = 0 = J(\alpha QX + a,\beta RY + b)$. Assume $\mathcal{A}\mathcal{B}>0$. By the homogeneity and isometry invariance of distance covariance \cite[Theorem~3(ii)]{Szekely2009} it follows that
\begin{align*}
    J(\alpha QX + a, \beta RY + b) &= \frac{2\operatorname{dCov}(\alpha QX+a,\beta RY+b)}{\sqrt{\operatorname{diam}(\operatorname{supp}(\alpha QX+a))\operatorname{diam}(\operatorname{supp}(\beta RY+b)})} \\
    &= \frac{2 \sqrt{\alpha\beta}\operatorname{dCov}(X,Y)}{\sqrt{\alpha\beta\mathcal{A}\mathcal{B}}} \\
    &= \frac{2\operatorname{dCov}(X, Y)}{\sqrt{\mathcal{A}\mathcal{B}}} \\
    &= J(X,Y).
\end{align*}

\bigskip
\noindent \textbf{(vi).} Let $\mathcal{A},\mathcal{B},\mathcal{C},\mathcal{D} \geq 0$ be the support diameters of $X,Y,Z,W$ respectively. Set
$\mathcal{S}:=\operatorname{diam}(\operatorname{supp}(X+Z))$ and
$\mathcal{T}:=\operatorname{diam}(\operatorname{supp}(Y+W))$.
Since $(X,Y) \ind (Z,W)$ we have $X \ind Z$ and $Y\ind W$. A standard support argument shows that
$\operatorname{supp}(X)+\operatorname{supp}(Z)\subseteq\operatorname{supp}(X+Z)$
and
$\operatorname{supp}(Y)+\operatorname{supp}(W)\subseteq\operatorname{supp}(Y+W)$. Hence
\begin{equation}
\label{ineq:sum-support-diameters}
    \mathcal{S}\geq\max\{\mathcal{A},\mathcal{C}\}, \qquad \mathcal{T}\geq\max\{\mathcal{B},\mathcal{D}\}.
\end{equation}
If $\mathcal{S}\mathcal{T} = 0$ then $J(X+Z,Y+W) = 0 \leq J(X,Y)+J(Z,W)$. Assume $\mathcal{S}\mathcal{T}>0$. Using the subadditivity of distance covariance under sums of independent pairs \cite[Theorem~3(iii)]{Szekely2009} and \eqref{ineq:sum-support-diameters} we get
{\allowdisplaybreaks
\begin{align*}
    J(X+Z,Y+W)
    &= \frac{2 \operatorname{dCov}(X+Z,Y+W)}{\sqrt{\mathcal{S}\mathcal{T}}} \\
    &\leq \frac{2}{\sqrt{\mathcal{S}\mathcal{T}}}
    (\operatorname{dCov}(X,Y)+\operatorname{dCov}(Z,W)) \\
    &= \sqrt{\frac{\mathcal{A}\mathcal{B}}{\mathcal{S}\mathcal{T}}}\,J(X,Y)
    + \sqrt{\frac{\mathcal{C}\mathcal{D}}{\mathcal{S}\mathcal{T}}}\,J(Z,W) \\
    &\leq J(X,Y) + J(Z,W).
\end{align*}
}
\end{proof}

It is known that the addition of independent noise to two random vectors cannot increase their distance covariance \cite{Szekely2009, Zhang2026}. A limitation of distance correlation as a normalization of distance covariance is that it does not have this property.

\begin{example}
Let $(X,Y)$ be distributed according to
\[
    {\renewcommand{\arraystretch}{1.2}
    \begin{array}{c|ccc}
        \hbox{\diagbox{Y}{X}} &           -1            &            1 &            2 \\ \hline \\[-2.75ex]
                          -2  & \phantom{-}0            & \frac{6}{16} & \frac{7}{16} \\
                \phantom{-}0  & \phantom{-}0            &            0 & \frac{1}{16} \\
                \phantom{-}1  & \phantom{-}0            & \frac{1}{16} &            0 \\
                \phantom{-}2  & \phantom{-}\frac{1}{16} &            0 &            0 \\[0.75ex]
                \hline
    \end{array}}
\]
\addvspace{0.75em}
\noindent Let $(Z,W) \ind (X,Y)$ with $Z \ind W$ and let $Z, W$ have marginal distributions
\[
    \Pr(Z=0)=\Pr(W=0)=\frac45,
    \qquad
    \Pr\!\left(Z=\frac1{10}\right)=\Pr\!\left(W=\frac1{10}\right)=\frac15.
\]
Applying \eqref{eq:dcov-identity} to all combinations of $X$ and $Y$ gives
\[
    \operatorname{dCov}^2(X,X)=\frac{1313}{4096},
    \qquad
    \operatorname{dCov}^2(Y,Y)=\frac{10865}{16384},
    \qquad
    \operatorname{dCov}^2(X,Y)=\frac{983}{8192},
\]
implying that
\[
    \operatorname{dCor}(X,Y)
    =
    \sqrt{\frac{983}{\sqrt{1313 \cdot 10865}}}
    \approx 0.5102.
\]
Similar calculations yield
\[
    \operatorname{dCov}^2(X+Z,X+Z)=\frac{19630881}{64000000},
    \qquad
    \operatorname{dCov}^2(Y+W,Y+W)=\frac{167487081}{256000000},
\]
\[
    \operatorname{dCov}^2(X+Z,Y+W)=\frac{15066337}{128000000},
\]
which imply that
\[
    \operatorname{dCor}(X+Z,Y+W)
    =
    \sqrt{\frac{15066337}{\sqrt{19630881 \cdot 167487081}}}
    \approx 0.5126.
\]
Thus for the tuple $(X,Y,Z,W)$
\[
    \operatorname{dCor}(X+Z,Y+W)>\operatorname{dCor}(X,Y),
\]
even though $Z$ and $W$ are independent of each other and of $(X,Y)$.
\end{example}

The following corollary of Proposition~\ref{prop:J-properties} assures that our normalization $J$ is nonincreasing under the addition of independent noise.
\begin{corollary}
\label{cor:J-independent-noise}
If $(X,Y)\ind (Z,W)$ and $Z\ind W$ then
\[
    J(X+Z,Y+W) \leq J(X,Y).
\]
\end{corollary}

\begin{proof}
By parts (vi) and (iii) of Proposition~\ref{prop:J-properties}
\[
    J(X+Z,Y+W) \leq J(X,Y) + J(Z,W) = J(X,Y).
\]
\end{proof}

\section{Sample bound and normalization}
\label{sec:sample-bound-normalization}
 
Distance covariance is commonly assessed through one of two estimators of its square: a biased V-statistic computed from double-centered distances \cite{Szekely2007} and an unbiased U-statistic computed from ``U-centered'' distances \cite{Szekely2014, Huo2016}.

Let $(x_1,y_1),\dots,(x_n,y_n)$ be a sample in $\mathbb{R}^N\times\mathbb{R}^M$. Define matrices $a$ and $b$ with entries $a_{ij}=\norm{x_i-x_j}$ and $b_{ij}=\norm{y_i-y_j}$ corresponding to raw pairwise distances. Writing the row, column, and grand means of the matrix $a$ respectively as
\[
    \bar{a}_{i\cdot} = \frac{1}{n}\sum_{j=1}^n a_{ij}, \qquad
    \bar{a}_{\cdot j} = \frac{1}{n}\sum_{i=1}^n a_{ij}, \qquad
    \bar{a}_{\cdot\cdot} = \frac{1}{n^2}\sum_{i,j=1}^n a_{ij},
\]
and similarly for $b$, the double-centered versions $\widehat{A}$ and $\widehat{B}$ of $a$ and $b$ are given by
\[
    \widehat{A}_{ij}=a_{ij}-\bar{a}_{i\cdot}-\bar{a}_{\cdot j}+\bar{a}_{\cdot\cdot},\qquad
    \widehat{B}_{ij}=b_{ij}-\bar{b}_{i\cdot}-\bar{b}_{\cdot j}+\bar{b}_{\cdot\cdot}.
\]
The usual sample distance covariance is defined as the nonnegative square root of the following V-statistic:
\begin{equation}
\label{eq:v-statistic}
    \operatorname{dCovV}^2_n(X,Y) = \frac{1}{n^2}\sum_{i,j=1}^n \widehat{A}_{ij}\widehat{B}_{ij}.
\end{equation}

Denote by $P_n=\frac{1}{n}\sum_{i=1}^n\delta_{(x_i,y_i)}$ the empirical measure of $(x_1,y_1),\dots,(x_n,y_n)$ and let $(X_n,Y_n)\sim P_n$. By Theorem~1 of \cite{Szekely2007}, $\operatorname{dCovV}^2_n(X,Y) = \operatorname{dCov}^2(X_n,Y_n)$. Define the sample diameters
\[
    \mathcal{A}_n := \max_{1 \leq i,j \leq n} a_{ij}, \qquad
    \mathcal{B}_n := \max_{1 \leq i,j \leq n} b_{ij}.
\]

\begin{corollary}
\label{cor:sample-bound}
$\operatorname{dCovV}_n(X,Y) \leq \frac{1}{2}\sqrt{\mathcal{A}_n \mathcal{B}_n}$.
\end{corollary}

\begin{proof}
By construction, $\operatorname{diam}(\operatorname{supp} X_n)=\mathcal{A}_n$ and
$\operatorname{diam}(\operatorname{supp} Y_n)=\mathcal{B}_n$. Applying Lemma~\ref{lem:sharp-bound} to $(X_n,Y_n)$ gives
\[
    \operatorname{dCov}(X_n,Y_n)\leq \frac{1}{2}\sqrt{\mathcal{A}_n\mathcal{B}_n}.
\]
The claim follows from $\operatorname{dCovV}_n(X,Y)=\operatorname{dCov}(X_n,Y_n)$.
\end{proof}

Corollary~\ref{cor:sample-bound} motivates the following estimator of the normalized measure of dependence defined in Section~\ref{sec:normalized-dependence-measure}:
\begin{equation}
\label{eq:sample-J-definition}
    J_n(X,Y) := \frac{2\operatorname{dCovV}_n(X, Y)}{\sqrt{\mathcal{A}_n\mathcal{B}_n}} = \frac{2\operatorname{dCov}(X_n, Y_n)}{\sqrt{\mathcal{A}_n\mathcal{B}_n}}, \quad \mathcal{A}_n \mathcal{B}_n > 0,
\end{equation}
with $J_n(X,Y):=0$ if $\mathcal{A}_n\mathcal{B}_n=0$, so that $J_n(X,Y) = J(X_n,Y_n)$. From this identity it is immediate that the empirical analogues of properties (i)--(v) in Proposition~\ref{prop:J-properties} hold for $J_n$. We show that, under i.i.d. sampling, $J_n(X,Y)$ is a strongly consistent estimator of $J(X,Y)$.

\begin{proposition}
\label{prop:J-consistency}
$J_n(X,Y) \to J(X,Y)$ a.s.
\end{proposition}

\begin{proof}
If $\mathcal{A}\mathcal{B}=0$ then $J(X,Y) = 0 = J_n(X,Y)$ a.s. for all $n$. Assume $\mathcal{A}\mathcal{B}>0$. 

By Theorem~2 of \cite{Szekely2007}, $\operatorname{dCovV}_n(X,Y)\to\operatorname{dCov}(X,Y)$ a.s. Since $\bigcup_{n\geq1}\operatorname{supp} X_n$ is dense in $\operatorname{supp} X$ a.s. and similarly for $Y$, we have
\[
    \mathcal{A}_n=\operatorname{diam}(\operatorname{supp} X_n)
    \to \operatorname{diam}(\operatorname{supp} X)=\mathcal{A}
    \quad \text{a.s.}
\]
and similarly for $\mathcal{B}_n$. Moreover since $\mathcal{A}\mathcal{B}>0$ we have with probability one that $\mathcal{A}_n\mathcal{B}_n>0$ for all sufficiently large $n$. By \eqref{eq:sample-J-definition}, \eqref{eq:J-definition}, and the continuity of the map $(u,a,b)\mapsto 2u/\sqrt{ab}$ on $[0,\infty)\times(0,\infty)^2$, it follows that
\[
    J_n(X,Y) = \frac{2\operatorname{dCovV}_n(X,Y)}{\sqrt{\mathcal{A}_n\mathcal{B}_n}} \to \frac{2\operatorname{dCov}(X,Y)}{\sqrt{\mathcal{A}\mathcal{B}}} = J(X,Y) \quad \text{a.s.}
\]
\end{proof}

We close by briefly considering the unbiased estimator of squared distance covariance introduced in \cite{Szekely2013a}. This quantity is often formulated as
\begin{equation}
\label{eq:u-statistic}
    \operatorname{dCovU}^2_n(X,Y)
    = \frac{1}{n(n-3)}\sum_{i\neq j}\widetilde{A}_{ij}\widetilde{B}_{ij}, \quad n \geq 4,
\end{equation}
where $\widetilde{A}, \widetilde{B}$ are the U-centered distance matrices \cite{Szekely2014}, obtained by adjusting the formulas for the entries of the double-centered matrices $\widehat{A}, \widehat{B}$ in the following way:
\begin{equation}
\label{eq:u-centered-distance}
    \widetilde{A}_{ij} = 
    \begin{cases}
        a_{ij} - \dfrac{n}{n-2}\bar{a}_{i\cdot} - \dfrac{n}{n-2}\bar{a}_{\cdot j} + \dfrac{n^2}{(n-1)(n-2)}\bar{a}_{\cdot\cdot} & i\neq j \\
        0 & i = j
    \end{cases}
\end{equation}
and similarly for $\widetilde{B}_{ij}$. Later shown to be a fourth-order U-statistic \cite{Huo2016}, this estimator may take negative values and therefore lack a real square root.

A practical workaround is to take the positive part of the estimator \cite{Geerligs2016, MonroyCastillo2025}. We show with a simple counterexample that the square root of this truncated statistic does not satisfy an inequality of the kind in Lemma~\ref{lem:sharp-bound}.

\begin{proposition}
The inequality
\[
\sqrt{\max\{\operatorname{dCovU}^2_n(X,Y),0\}}
\leq\frac{1}{2}\sqrt{\mathcal{A}_n\mathcal{B}_n}
\]
does not in general hold.
\end{proposition}

\begin{proof}
Consider the sample
\[
    x = (0,0,1,1), \qquad y = (0,0,1,1).
\]
The raw distance matrices $a, b$ satisfy $a_{ij} = b_{ij} = 0$ for all index pairs lying within either set $\{1,2\}$ or $\{3,4\}$ and $a_{ij} = b_{ij} = 1$ for all index pairs lying between those sets. Hence every row and column of $a$ and $b$ contains two zeros and two ones, so that their row, column, and grand means all equal $1/2$. Using \eqref{eq:u-centered-distance} gives $\widetilde{A}_{ij} = \widetilde{B}_{ij} = -2/3$ for the four off-diagonal index pairs lying within the two sets and $\widetilde{A}_{ij} = \widetilde{B}_{ij} = 1/3$ for the eight off-diagonal index pairs lying between them. Therefore $\sum_{i\neq j}\widetilde{A}_{ij}\widetilde{B}_{ij} = 8/3$ and so
$\operatorname{dCovU}^2_{4}(X,Y) = 2/3$. Here $\mathcal{A}_4 = \mathcal{B}_4 = 1$, and thus
\[
    \sqrt{\max\{\operatorname{dCovU}^2_4(X,Y),0\}} = \sqrt{\frac{2}{3}} > \frac{1}{2} = \frac{1}{2}\sqrt{\mathcal{A}_4\mathcal{B}_4}.
\]
\end{proof}

\section{Concluding remarks}
\label{sec:concluding-remarks}

We obtain a sharp upper bound on distance covariance that depends only on the diameters of the supports of the underlying random vectors. The bound obtained gives the maximum distance covariance under a pair of support diameter constraints. Dividing distance covariance by this maximum yields a normalized measure of dependence that combines the scale invariance of distance correlation with the subadditivity of distance covariance under sums of independent pairs.

While the sum-subadditivity property in part (vi) of Proposition~\ref{prop:J-properties} and the independent-noise property in Corollary~\ref{cor:J-independent-noise} rely on the linear structure of Euclidean space, the proof of the upper bound given in Lemma~\ref{lem:sharp-bound} uses only the metric structure of the random vector supports. Hence the sharp bound and the related normalization of distance covariance extend directly to bounded random elements in separable metric spaces of strong negative type \cite{Lyons2013,Lyons2018}.

\bibliographystyle{amsalpha}
\bibliography{main}

\end{document}